\documentclass[11pt,reqno]{amsart}

\usepackage[PostScript=dvips,dpi=600,nohug,small]{diagrams}
\usepackage{amssymb}
\usepackage{mathrsfs}
\usepackage{cite}
\usepackage{amsfonts}
\usepackage[symbol*]{footmisc}
\newcommand{\rmv}[1]{}

\def\wh{\widehat}

\def\pv#1{\ensuremath{{\bf#1}}}

\def\ilim{\varprojlim}

\def\inv{^{-1}}
\def\p{\varphi}

\def\pinv{{\p \inv}}
\def\J{\mathrel{{\mathscr J}}} % J - relation
\def\R{\mathrel{{\mathscr R}}} % R - relation
\def\L{\mathrel{{\mathscr L}}} % L - relation
\def\e<{\leq _{E}}

\def\ov#1{\ensuremath{\overline {#1}}}

\def\til#1{\ensuremath{\widetilde {#1}}}

\def\malce{\mathbin{\hbox{$\bigcirc$\rlap{\kern-8.3pt\raise0,50pt\hbox{$\mathtt{m}$}}}}}
\def\FP#1#2{\ensuremath{\widehat{F_{\mathbf #1}}(#2)}}

\def\1sk{^{(1)}}

\def\to{\rightarrow}

\def\hatexp#1#2{\widetilde{#2}}

\def\Thmname{Theorem}
\def\Propname{Proposition}
\def\Lemmaname{Lemma}
\def\Definitionname{Definition}
\newtheorem{Thm}{\Thmname}%[section]

\newtheorem{Lemma}[Thm]{\Lemmaname}

\newtheorem{Cor}[Thm]{Corollary}

\numberwithin{equation}{section}

\title[Closed subgroups of free profinite monoids]{Closed subgroups of free profinite monoids are projective profinite groups}

\author{John Rhodes\and Benjamin Steinberg}
\address{Department of Mathematics\\
University of California at Berkeley \\
Berkeley \\ CA 94720\\
USA\\ \and School of Mathematics and Statistics\\
Carleton University \\
1125 Colonel By Drive\\
Ottawa, Ontario  K1S 5B6 \\
Canada}
\thanks{The second author was supported in part by NSERC}
\email{rhodes@math.berkely.edu\and bsteinbg@math.carleton.ca}
\date{\today}

\begin{document}
\begin{abstract}
We prove that the class of closed subgroups of free profinite
monoids is precisely the class of projective profinite groups. In
particular, the profinite groups associated to minimal symbolic
dynamical systems by Almeida are projective. Our result answers a question
raised by Lubotzky during the lecture of Almeida at the Fields
Workshop on Profinite Groups
and Applications, Carleton University, August~2005.  We
also prove that any finite subsemigroup of a free profinite monoid
consists of idempotents.
\end{abstract}
\maketitle

\section{Introduction}

It has long been an open question whether closed subgroups of free
profinite monoids must be projective profinite groups. At one time
it was even hoped that maximal subgroups would be free profinite
groups.
  Almeida recently
associated to each minimal symbolic dynamical system a maximal
subgroup of a free profinite monoid, which is a conjugacy invariant
of the system~\cite{Jorgesubgroup}.  Almeida gave sufficient
conditions for the associated subgroup to be a free profinite group
and, in particular, showed this to be the case for the group
associated to a Sturmian system~\cite{Jorgesubgroup}.  He also
provided an example of a dynamical system leading to a maximal
subgroup that is not free profinite~\cite{Jorgesubgroup}.  Almeida
lectured about this work at the Fields Workshop on Profinite Groups
and Applications, Carleton University, August~2005.  Lubotzky asked
at the end of the talk whether all closed subgroups of a free
profinite monoid are projective. In this paper, we provide a
positive answer to Lubotzky's question by showing that closed
subgroups of free profinite monoids are indeed projective. In fact,
a more general result holds.

Almeida and Weil~\cite{AWsurvey,MR1297147} defined a pseudovariety
of groups $\pv H$ to be \emph{arborescent} if whenever $1\to A\to
G\to H\to 1$ is a short exact sequence of groups with $A,H\in \pv H$
and $A$ abelian, one has $G\in \pv H$. The origin of the terminology
is that if $\pv H$ is non-trivial, then being arborescent is
equivalent to having the Cayley graph of each free pro-\pv H group
being a pro-\pv H tree~\cite{AWsurvey,MR1297147}.

If \pv H is a pseudovariety of groups, then $\ov {\pv H}$ denotes
the pseudovariety of monoids whose subgroups belong to $\pv H$. This
note proves that if $\pv H$ is an arborescent pseudovariety of
groups, then any closed subgroup of a free pro-$\ov {\pv H}$ monoid
is a projective profinite group and that any finite subsemigroup of
a free pro-$\ov {\pv H}$ monoid is an idempotent semigroup.  In
particular, by considering the pseudovariety of all finite groups
and the trivial pseudovariety, this result applies to free profinite
monoids and free profinite aperiodic monoids. A semigroup is termed
\emph{aperiodic} if all its subgroups are trivial.

Before stating precise results, let us first recall some
definitions. A profinite group $G$ is called
\emph{projective}~\cite{RZbook} if whenever one has a diagram of groups
(called an \emph{embedding problem}~\cite{RZbook})
\begin{equation}\label{embeddingproblem}
\begin{diagram}
 &      &   &      &  &              &  G     &  &\\
 &      &   &      &  &              &  \dOnto_{\p}&     &\\
1& \rTo & K & \rTo & A& \rTo^{\alpha}& B  &\rTo & 1
\end{diagram}
\end{equation}
with $A$ finite and $\p$ a continuous epimorphism, there is a
continuous lift $\til \p:G\to A$, called a \emph{weak solution} to
the embedding problem~\cite{RZbook}, such that
\begin{equation}\label{findlift}
\begin{diagram}
 &      &   &      &  &              &  G     &  &\\
 &      &   &      &  &  \ldTo^{\til \p}            &  \dOnto_{\p}&     &\\
1& \rTo & K & \rTo & A& \rTo^{\alpha}& B  &\rTo & 1
\end{diagram}
\end{equation}
 commutes. If \pv H is a pseudovariety of
groups, then a pro-\pv H group is called \pv H-projective if each
embedding problem \eqref{embeddingproblem} with $A\in \pv H$ has a
weak solution \eqref{findlift}; see~\cite{RZbook}.  It turns out
that for arborescent pseudovarieties \pv H, being projective and \pv
H-projective coincide (see Corollary~\ref{equivalence}).

If $\pv V$ and $\pv W$ are pseudovarieties of monoids, then  $\pv
V\ast \pv W$ denotes the pseudovariety of monoids generated by
semidirect products of the form \mbox{$V\rtimes W$} with $V\in \pv
V$ and $W\in \pv W$. If $\pv V$ is a pseudovariety of semigroups and
$\pv W$ is a pseudovariety of monoids, then their Mal'cev product
$\pv V\malce \pv W$ is the pseudovariety generated by monoids $M$
admitting a homomorphism  $\p:M\to N$ with $N\in \pv W$ such that
$e\pinv\in \pv V$ for each idempotent $e\in N$.  A finite semigroup
$S$ is called \emph{completely regular} if each element $s\in S$
satisfies $s^m=s$ for some $m>0$.  A finite semigroup is called a
\emph{band} if each element is idempotent.  A band is precisely the
same thing as an aperiodic completely regular semigroup.  Let $\pv
{Ab}$ denote the pseudovariety of finite abelian groups and \pv A
denote the pseudovariety of finite aperiodic monoids.  Our main
results are then the following two theorems.

\begin{Thm}\label{theorem2}
Let $\pv V$ be a pseudovariety of monoids with the property that
$(\pv V\cap \pv {Ab})\ast \pv V = \pv V$. Then every closed subgroup
of a free pro-\pv V monoid is a projective profinite group, and
in particular torsion-free.
\end{Thm}

\begin{Thm}\label{theorem1}
Let $\pv V$ be a pseudovariety of monoids such that $\pv A\malce \pv
V=\pv V$.  Then every finite subsemigroup of a free pro-\pv V monoid
is completely regular.
\end{Thm}

Since pseudovarieties of the form $\ov {\pv H}$ with $\pv H$
arborescent evidently satisfy the hypotheses of
Theorems~\ref{theorem2} and~\ref{theorem1}, we deduce our main
result:

\begin{Cor}\label{mainresult}
If $\pv H$ is an arborescent pseudovariety of groups, then every
closed subgroup of a free pro-$\ov{\pv H}$ monoid is a projective
profinite group and each finite subsemigroup of a free pro-$\ov{\pv
H}$ monoid is a band.
\end{Cor}

It is well known that a profinite group is projective if and only if
it is a closed subgroup of a free profinite group~\cite{RZbook}.  As
the natural continuous projection from the free profinite monoid
onto the free profinite group
splits~\cite{AlmeidaVolkov,RhodesStein}, it in fact follows that
projective profinite groups are \emph{precisely} the closed
subgroups of free profinite monoids.

The paper is organized as follows.  The first section gives a proof
that closed subgroups of
projective groups are projective, highlighting the main idea
underlying Theorem~\ref{theorem2}.  The next two sections of the
paper prove Theorems~\ref{theorem2} and~\ref{theorem1}.  The final
section gives an application to projective semigroups.

\section{Projective profinite groups}
The usual proof (c.f.~\cite{RZbook}) that a closed subgroup of a
projective profinite group is projective goes roughly as follows.
First one shows that a projective profinite group embeds in a free
profinite group.  Then one shows that projectivity of closed
subgroups of free profinite groups reduces to the case of clopen
subgroups. Clopen subgroups are then shown to be free using Hall's
theorem~\cite{Hall} and the Nielsen-Schreier theorem.  We discovered
an elementary proof of this result using the monomial
map~\cite{Hallbook,Weilandt} (sometimes called the
Krasner-Kaloujnine embedding~\cite{KrasnerKal}). Luis
Ribes\footnote{We are grateful to Luis Ribes for directing our
attention to~\cite{Cossey}.} pointed out to us that this same proof
scheme was used in~\cite{Cossey}. Nonetheless, we provide the proof
as it relies on two lemmas that we also use in the monoidal context.
The first lemma contains the \emph{key idea} for \emph{both} the
group and monoid cases.

\begin{Lemma}\label{thetrick}
Let $G$ be a profinite group and suppose that one has an embedding
problem \eqref{embeddingproblem} which has no weak solution. Suppose
that $\til B$ is a finite group with $B\leq \til B$.  Then one can
find a diagram of finite groups
\begin{equation}\label{embeddingproblem2}
\begin{diagram}
 &      &   &      &  &              &  G     &  &\\
 &      &   &      &  &              &  \dTo_{\p}&     &\\
1& \rTo & \til K & \rTo & \til A& \rTo^{\til{\alpha}}& \til B  &\rTo
& 1
\end{diagram}
\end{equation}
such that $\ker \til{\alpha}$ is a subgroup of $K^n$ for some $n$
and there is no continuous homomorphism $\til {\p}:G\to \til A$
lifting $\p$ (i.e.\ so that $\til {\p}\til{\alpha} = \p$).
\end{Lemma}
\begin{proof}
Let $H$ be the quotient of $\til B$ by the kernel of the action on
the set $\til B/B$ of right cosets of $B$.  Then there is a
well-known embedding
\[\til B\hookrightarrow B\wr (\til B/B,H) = B^{\til B/B}\rtimes
H,\] see~\cite{Hallbook,Weilandt,Eilenberg,qtheor,KrasnerKal}.  The
construction is as follows.  Choose coset representatives $\ov {Bg}$
for each right coset of $B$ so that $\ov B=1$.  Let us write $[g]$
for the image in $H$ of $g\in \til B$.  Then the embedding takes
$g\in \til B$ to $(f_g,[g])$ where $(Bb)f_g = \ov {Bb}g(\ov
{Bbg})\inv$.  In particular, if $g\in B$, then $Bf_g = g$.  Thus the
map $(f_g,[g])\mapsto Bf_g$ is an isomorphism $\tau$ from the copy
of $B$ in $\til B$ to $B$.

There is a natural epimorphism $\psi:A^{\til B/B}\rtimes
H\twoheadrightarrow B^{\til B/B}\rtimes H$ induced by $\alpha$,
namely $(h,[g])\psi = (h\alpha,[g])$.  Also $\ker \psi \cong K^{\til
B/B}$ since it consists of all pairs $(h,[1])$ with $(Bg)h\in \ker
\alpha=K$, all $Bg\in \til B/B$. So let $\til A = \til B\psi\inv$
and $\til{\alpha} = \psi|_{\til A}$. Then $\ker\til{\alpha}$ is
isomorphic to a subgroup of $K^{\til B/B}$. Let $A' =
B\til{\alpha}\inv$. We claim that the map $\rho:A'\to A$ given by
$(f,[g])\mapsto Bf$ is a homomorphism such that the diagram
\begin{equation}\label{needcommmute}
\begin{diagram}
A'& \rTo^{\til{\alpha}}& B &\rInto & B^{\til B/B}\rtimes H\\
\dTo^{\rho} & & \dTo_{\tau} \\
A & \rTo^{\alpha} & B
\end{diagram}
\end{equation}
commutes.  Let us assume this for the moment and complete the proof.
If $\til{\p}$ lifts $\p$ in \eqref{embeddingproblem2}, then
$\til{\p}\rho$ solves the original embedding problem
\eqref{embeddingproblem}.  Therefore, $\p$ has no lift by the
hypothesis.

Clearly \eqref{needcommmute} commutes by definition of $\rho$ and
$\tau$, so we just need to check that $\rho$ is a homomorphism.
Indeed, if $(f,[g])$ and $(f',[g'])$ are in $A'$ (and so we may take
$g,g'\in B$), then $(f,[g])(f',[g']) = (f{}^{[g]}\!{f'},[gg'])$. As
$g\in B$,
\[((f,[g])(f',[g']))\rho = B(f{}^{[g]}\!{f'}) = Bf\cdot Bgf' = Bf\cdot
Bf' = (f,[g])\rho (f',[g'])\rho.\]  This completes the proof.
\end{proof}

Our next lemma is a reduction on the types of embedding problems one
must solve to establish projectivity.

\begin{Lemma}\label{gocofinal}
Let $G=\ilim _{i\in I} G_i$ with the $G_i$ finite continuous
quotients of $G$ and let $\pi_i:G\twoheadrightarrow G_i$ be the
projection. Consider an embedding problem as per
\eqref{embeddingproblem}.  Then there is an index $i$  and an
epimorphism $\psi:G_i\to B$ such that the diagram
\begin{equation}\label{cofinaldiag}
\begin{diagram}
                   &                 &              &                & &G\\
                   &                 &              & &\ldOnto(3,2)^{\pi_i}\ldOnto(3,4)_{\p}  &\\
A\times _{\psi} G_i&\rOnto^{\alpha'} & G_i          &  &              &\\
\dOnto^{\psi^*}             &                 & \dOnto_{\psi}&   &             & \\
A                   &\rOnto^{\alpha}  & B            &&&
\end{diagram}
%\begin{diagram}
% A\times _{\psi} G_i&\rOnto^{\alpha'} & G_i           & \lTo^{\pi_i}   & G\\
% \dOnto^{\psi^*}     &               & \dOnto_{\psi } & \ldTo_{\p}    & \\
% A                 &\rOnto^{\alpha}  & B            &    &
% \end{diagram}
\end{equation}
commutes, where $\psi^*:A\times _{\psi} G_i\to A$ is the pullback of
$\psi$ along $\alpha$.  Moreover, $\ker \alpha'\cong \ker \alpha$.

In particular, $G$ is a projective profinite group if and only if
all embedding problems \eqref{embeddingproblem} with $B=G_i$ and
$\p=\pi_i$, some $i\in I$, have a weak solution.
\end{Lemma}
\begin{proof}
Since $B$ is finite, $\p$ factors through one of the continuous
projections $\pi_i:G\twoheadrightarrow G_i$ via an epimorphism
$\psi:G_i\twoheadrightarrow B$. Consider the pullback $A\times
_{\psi} G_i= \{(a,g)\in A\times G_i\mid a\alpha = g\psi\}$ of $\psi$
along $\alpha$.  Then the projections $\alpha'$ to $G_i$ and
$\psi^*$ to $A$ are epimorphisms as $\alpha$ and $\psi$ are
epimorphisms. Also, $\ker \alpha' = \{(a,1)\mid a\alpha = 1\psi
=1\}\cong \ker \alpha$. This proves the existence of
\eqref{cofinaldiag}.

The final statement follows as a weak solution to the embedding
problem for $\alpha'$ in \eqref{cofinaldiag} leads to a weak
solution to the embedding problem \eqref{embeddingproblem}.
\end{proof}

\begin{Cor}\label{closedundersubgroups}
Closed subgroups of projective profinite groups are projective.  In
particular, projective profinite groups are torsion-free.
\end{Cor}
\begin{proof}
Suppose that $H$ is a closed subgroup of a projective profinite
group $G$.  Write $G = \ilim_{i\in I} G_i$ with the $G_i$ finite
quotients and let $\pi_i:G\twoheadrightarrow G_i$ be the projection.  As $H$ is
closed,  $H = \ilim _{i\in I} H_i$ where $H_i = H\pi_i$.  Suppose
that $H$ is not projective. Then by Lemma~\ref{gocofinal}, there is
an embedding problem
\[\begin{diagram}
 &      &   &      &  &              &  H     &  &\\
 &      &   &      &  &              &  \dOnto_{\pi_i}&     &\\
1& \rTo & K & \rTo & A& \rTo^{\alpha}& H_i  &\rTo & 1
\end{diagram}\]
having no weak solution.  Taking $B= H_i$ and $\til B =G_i$ in
Lemma~\ref{thetrick}, we obtain a diagram
\[\begin{diagram}
 &      &   &      &  &              &  H     &  &\\
 &      &   &      &  &              &  \dTo_{\pi_i}&     &\\
1& \rTo & \til K & \rTo & \til A& \rTo^{\til{\alpha}}& G_i  &\rTo &
1
\end{diagram}\]
so that $\pi_i$ cannot be lifted to $\til A$.  It follows that the
embedding problem
\[\begin{diagram}
 &      &   &      &  &              &  G     &  &\\
 &      &   &      &  &              &  \dOnto_{\pi_i}&     &\\
1& \rTo & \til K & \rTo & \til A& \rTo^{\til{\alpha}}& G_i  &\rTo &
1
\end{diagram}\]
has no weak solution, contradicting that $G$ is projective.

Since any finite subgroup of a projective profinite group is closed,
it follows that such a subgroup must be projective.   Since a cyclic
group of prime order $p$ is not projective (as the canonical
projection $\mathbb Z/p^2\mathbb Z\to \mathbb Z/p\mathbb Z$ does not
split) we conclude that projective profinite groups are
torsion-free.
\end{proof}

Now we wish to show that if \pv H is an arborescent pseudovariety of
groups, then the \pv H-projective groups are projective and that one
only needs to consider very special embedding problems. According
to~\cite{RZbook}, a pseudovariety of groups \pv H  is called
\emph{saturated} if whenever $G$ is a finite group such that its
Frattini quotient $G/\Phi(G)$ belongs to $\pv H$, then $G\in \pv H$.
It is shown in~\cite[Satz III.3.8]{Huppert} (see also~\cite[Example
7.6.5]{RZbook}) that if $p$ is a prime dividing the order of the
Frattini subgroup $\Phi(G)$, then $p$ divides the order of the
Frattini quotient $G/\Phi(G)$. Since $\Phi(G)$ is nilpotent it
follows that if \pv H is arborescent, then $\pv H$ is saturated.

The importance of saturated pseudovarieties is that groups from
these pseudovarieties lift.  The following is~\cite[Lemma
7.6.6]{RZbook}.

\begin{Lemma}\label{saturated}
Let \pv H be a saturated pseudovariety of groups and let $\p:G\to H$
be an epimorphism of finite groups with $H\in \pv H$.  Then there
exists $M\leq G$ such that $M\in \pv H$ and $M\p = H$.
\end{Lemma}

The next corollary follows from well-known results, which can be
found in~\cite{RZbook}, that will be used for our monoidal results.

\begin{Cor}\label{equivalence}
Let $\pv H$ be an arborescent pseudovariety of groups and let $G$ be
a pro-\pv H group.  Then the following are equivalent:
\begin{enumerate}
\item $G$ is projective;
\item $G$ is \pv H-projective;
\item All embedding problems \eqref{embeddingproblem} for $G$ with
$K,B\in \pv H$ and $K$ an elementary abelian $p$-group have a weak
solution.
\end{enumerate}
\end{Cor}
\begin{proof}
The implications (1) implies (2) implies (3) are trivial (the last
uses that \pv H is arborescent and so $A\in \pv H$).  For (3)
implies (1), we use that in order to show that $G$ is projective, it
suffices by~\cite[Theorem 7.5.1]{RZbook} and~\cite[Proposition
7.5.4]{RZbook} to establish the existence a weak solution for all
embedding problems \eqref{embeddingproblem} for $G$ with $A$ finite
and $K$ an elementary abelian $p$-group. But given such an embedding
problem, we have $B=G\p\in \pv H$ since $G$ is pro-\pv H.
Lemma~\ref{saturated} allows us to replace $A$ by a subgroup $M$
belonging to \pv H.  Then $M\cap K$ is an elementary abelian
$p$-group in \pv H, as required.
\end{proof}

\section{Proof of Theorem ~\ref{theorem2}}
We proceed to prove Theorem~\ref{theorem2}.  If $X$ is a topological
space, then we denote by $\FP V X$ the free pro-\pv V monoid on
$X$~\cite{Almeida:book,AWsurvey,qtheor}.  Let $G$ be a closed
subgroup of $\FP V X$.  Denote by $\pv H$ the pseudovariety of all
groups belonging to \pv V. Then \pv H is arborescent and $G$ is
pro-\pv H. Suppose that $G$ is not projective. Then by
Corollary~\ref{equivalence}, there is an embedding
problem~\eqref{embeddingproblem}, with $K,B\in \pv H$ and $K$ an
elementary abelian $p$-group, having no weak solution. Since $G$ is
closed, $G$ is the projective limit of its images under the inverse
system of continuous finite quotients of $\FP V X$.
Lemma~\ref{gocofinal} then implies that we may assume without loss
of generality that $\p:G\to B$ is the restriction of a continuous
surjective homomorphism $\FP V X\twoheadrightarrow M$ with $M\in \pv
V$.

Let $e\in B$ be the identity and let $R=\{m\in M\mid mM=eM\}$ denote
the $\R$-class~\cite{Green} of
$B$ in $M$. We consider the Sch\"utzenberger representation of
$M$ on $R$. Namely, $M$ acts on the right of $R$ by \emph{partial}
transformations via
\begin{equation}\label{schutzrepcoordfree}
r\cdot m =
\begin{cases} rm &
  rm\in R\\ \text{undefined} &rm\notin R\end{cases}
\end{equation}
    Let $N$ be
the quotient of $M$ by the kernel of this action.  Then $B$ maps
injectively into $N$, since it clearly acts faithfully on itself and
hence on $R$. Therefore, by replacing $M$ with $N$,  we may assume
without loss of generality that $M$ acts faithfully on the right of
$R$.

Let $H$ be the maximal subgroup of $M$ containing $B$.  So $H$ is
the group of units of the monoid $eMe$.  By Lemma~\ref{thetrick}
with $\til B = H$, we can find a diagram
\begin{equation}\label{almostdonewithproj}
\begin{diagram}
 &      &   &      &  &              &  G     &  &\\
 &      &   &      &  &              &  \dTo_{\p}&     &\\
1& \rTo & \til K & \rTo & \til A& \rTo^{\til{\alpha}}& H  &\rTo & 1
\end{diagram}
\end{equation}
such that $\p$ has no lift and $\til K$ is an elementary abelian
$p$-group (being a subgroup of a direct power of $K$).

In order to reuse the \emph{key idea} from the proof of
Lemma~\ref{thetrick}, we need an embedding of $M$ into a wreath
product of $H$ with some monoid.  The embedding we use is due to
Sch\"utzenberger~\cite{Schutzrep,Schutzmonomial}.
In~\cite{Schutzrep}, Sch\"utzenberger proved that $H$ acts freely on
the left of $R$ by automorphisms of the right action
\eqref{schutzrepcoordfree} of $M$ on $R$ ; see~\cite[Theorem
2.22]{CP} and~\cite{qtheor}. Hence $M$ acts on the right of the
orbit set $Q= H\backslash R$ by partial
transformations.\footnote{The set $Q$ is in bijection with the set
of $\L$-classes of the $\J$-class of
$H$~\cite{Green,CP,Arbib,qtheor}.} Set $Q^{\star} = Q\cup \{\star\}$
where $\star\notin Q$. One can define an action of $M$ on the right
of $Q^{\star}$ by total functions by defining, for $r\in R$ and
$m\in M$,
 \[Hr\cdot m = \begin{cases} Hrm & rm\in R\\ \star & rm\notin
 R\end{cases}\] and by defining $\star\cdot m=\star$.  Let $N$ be the
 quotient of $M$ by the kernel of this action.  Let $H^0=H\cup \{0\}$
(viewed, say, as a subsemigroup of the group ring of $H$).  Then
there is an embedding $M\hookrightarrow H^0\wr (Q^{\star},N) =
(H^0)^{Q^\star}\rtimes N$, where $n\in N$ acts on $f:Q^{\star}\to
H^0$ by $x{}^{n}\!{f} = xnf$ (see~\cite[Proposition 8.2.17]{Arbib}).

Let us briefly describe the embedding.  We use $[m]$ for the image
of $m\in M$ in $N$.  Choose a representative $\ov {Hr}$ for each
orbit of $H$ on $R$ so that $\ov H = e$.   Then $m\mapsto (f_m,[m])$
where $\star f_m = 0$ and, for $r\in R$, $Hrf_m = 0$ if $Hr\cdot m=\star$
and otherwise is the unique element of $H$ such that $\ov {Hr}\cdot m =
(Hr)f_m \ov {Hr\cdot m}$.  The existence and uniqueness of this element comes from the
fact that $H$ acts freely on the left of $R$ by automorphisms of the
right action \eqref{schutzrepcoordfree}. Notice that
if $h\in H$, then $h\mapsto (f_h,[h])$ where $Hf_h=h$ since $\ov
H\cdot h = h = h\ov H = h\ov {H\cdot h}$.  That is, the map $\tau$
sending $(f_h,[h])$ to
$Hf_h$ is an isomorphism from the copy of $H$ in $M\leq H^0\wr
(Q^\star,N)$ to $H$.

  Let $\til{\alpha}^0:\til
A^0\twoheadrightarrow H^0$ be the induced surjective homomorphism.
The reader easily verifies that each local monoid of the derived
category\footnote{The reader just interested in the proof of
Corollary~\ref{mainresult} may ignore arguments involving derived
categories.}  $D_{\til{\alpha}^0}$ of $\til{\alpha}^0$ is $\til
K=\ker \til{\alpha}$ except the local monoid at $0$, which is
trivial. Hence $D_{\til{\alpha}^0}$ divides an elementary abelian
$p$-group by the locality of group pseudovarieties~\cite{Tilson}.

Next, we consider the natural surjective homomorphism
\[\psi:\til A^0\wr (Q^\star,N)\twoheadrightarrow H^0\wr (Q^\star,N)\] induced
by $\til{\alpha}^0$ ---  so $(f,n)\mapsto (f\til{\alpha}^0,n)$. Set
$M' = M\psi\inv$ and consider the surjective homomorphism
$\lambda=\psi|_M:M'\twoheadrightarrow M$.  Since
$D_{\til{\alpha}^0}$ divides an elementary abelian
$p$-group,~\cite[Theorem 7.1]{Cats2} implies that the derived
category of $\psi$, and hence of its restriction $\lambda$
(using~\cite[Proposition 5.12]{Cats2}), divides an elementary
abelian $p$-group.  As $\pv V$ contains a cyclic group of order $p$
and $M\in \pv V$, we have $M'\in (\pv V\cap\pv {Ab})\ast \pv V=\pv
V$ by the Derived Category Theorem~\cite{Tilson}\footnote{If one is
just interested in Corollary~\ref{mainresult}, then it is easy to
verify directly that each subgroup of $M'$ is an extension of an
elementary abelian $p$-group by a subgroup of $M$.}.

Let $A'=H\lambda\inv$ be the preimage of $H$ in $M'$.   We claim
that the map $\rho:A'\to \til A$ given by $(f,[h])\mapsto Hf$ is a
homomorphism such that the diagram
\begin{equation}\label{needcommmute2}
\begin{diagram}
A'& \rTo^{\lambda}& H &\rInto & H^0\wr (Q^\star,N)\\
\dTo^{\rho} & & \dTo_{\tau} \\
\til A & \rTo^{\til{\alpha}} & H
\end{diagram}
\end{equation}
commutes.  Let us assume this for the moment and complete the proof.
Since $\FP V X$ is free, the natural projection $\FP V
X\twoheadrightarrow M$ (whose restriction to $G$ is $\p$) can be
lifted to a continuous homomorphism $\eta:\FP V X\to M'$. Then
$G\eta\leq A'$ and $\eta\rho$ lifts $\p$, contradicting the
non-existence of a lift in \eqref{almostdonewithproj}.

Since \eqref{needcommmute2} commutes by definition of $\rho$ and
$\tau$, we just need to check that $\rho$ is a homomorphism. Well,
if $(f,[h])$ and $(f',[h'])$ are in $A'$ (and so we may take
$h,h'\in H$), then $(f,[h])(f',[h']) = (f{}^{[h]}\!{f'},[hh'])$.
Since $h\in H$, we have \[((f,[h])(f',[h']))\rho =
H(f{}^{[h]}\!{f'}) = Hf\cdot Hhf' = Hf\cdot Hf' = (f,[h])\rho
(f',[h'])\rho.\]  This completes the proof of
Theorem~\ref{theorem2}.\qed

\section{Proof of Theorem~\ref{theorem1}}
This proof uses the Henckell-Sch\"utzenberger expansion
(see~\cite{BR--exp,ourstablepairs} for details).   If $X$ is a set,
then $X^*$ denotes the free monoid on $X$.  If $M$ is an
$X$-generated monoid, then one can define a congruence $\equiv$ on
$X^*$ by $w_1\equiv w_2$ if and only if, for each factorization
$w_1=u_1v_1$, there is a factorization $w_2=u_2v_2$ with
$[u_1]_M=[u_2]_M$, $[v_1]_M=[v_2]_M$, and conversely.  The quotient
$\hatexp 2 M =M/\!{\equiv}$ is finite if $M$ is finite and there is
a natural surjective homomorphism $\eta:\hatexp 2
M\twoheadrightarrow M$ of $X$-generated monoids. Moreover, $\eta$ is
an aperiodic morphism, meaning that $e\eta\inv \in \pv A$ for each
idempotent $e\in M$ (see~\cite{BR--exp,ourstablepairs}).  In
particular, if $\pv A\malce \pv V=\pv V$, then $M\in \pv V$ implies
$\hatexp 2 M\in \pv V$.  The monoid $\hatexp 2 M$ is called the
\emph{Henckell-Sch\"utzenberger} expansion of $M$.

  Let $S$ be a
finite subsemigroup of $\FP V X$.  By finiteness, one can find a
continuous surjective homomorphism $\p:\FP V X\to M$ with $M\in \pv
V$ such that $\p|_S$ is injective.  If $Y=X\p$ (which is a finite
set), then $\p$ must factor through $\FP V Y$ and so we deduce that
$S$ is a finite subsemigroup of $\FP V Y$.  So without loss of
generality, we may assume that $X$ is finite.

Suppose that $S$ is not completely regular.   We shall arrive at a
contradiction.  As $S$ is not completely regular, there is an
element $\sigma\in S$, which is not a group element, such that
$\sigma^2$ is a group element.  Then $\sigma^2 = \sigma^{2+k}$ for
some $k>0$.  Without loss of generality we may assume that $k\geq
2$.   Now $\hatexp 2 M\in \pv V$ by our hypothesis and so $\p$ must
factor continuously through $\eta:\hatexp 2 M\to M$.

 Since $X^*$ is dense in $\FP V X$, we can find a word $w\in
X^*$ such that $[w]_{\hatexp 2 M} = [\sigma]_{\hatexp 2 M}$. Then we
have $[w^2]_{\hatexp 2 {M}}=[\sigma^2]_{\hatexp 2
{M}}=[\sigma^k]_{\hatexp 2 M}=[w^k]_{\hatexp 2 M}$. Hence there is a
factorization $w^k = w^{k_1}xyw^{k_2}$ with $x,y\in X^*$ such that
$w = xy$, $k_1,k_2\geq 0$, $k_1+k_2+1=k$ and $[w^{k_1}x]_{M} =
[w]_{M}$ and $[yw^{k_2}]_{M} = [w]_{M}$.  Therefore, we have the
equations:
\begin{equation}\label{completelyreg}
\begin{split}
[\sigma^{k_1}]_{M}[x]_{M} &= [\sigma]_{M}\\
[y]_{M}[\sigma^{k_2}]_{M} &= [\sigma]_{M}.
\end{split}
\end{equation}
As $k_1+k_2+1=k\geq 4$, it cannot be the case that both $k_1,k_2\leq
1$. We conclude from \eqref{completelyreg} that $[\sigma]_{M} \J
[\sigma^2]_{M}$ and so $[\sigma]_M$ is a group
element~\cite{Green,Arbib,CP}.  Since $\p|_S$ is injective, $\sigma$
itself is a group element, a contradiction. This completes the proof
of Theorem~\ref{theorem1}.\qed

\section{An application to projective semigroups}
A finite semigroup $S$ is said to be \emph{projective}~\cite{qtheor}
if any onto homomorphism $\p:T\twoheadrightarrow S$, with $T$
finite, splits (i.e.\ there exists $\psi:S\to T$ with
$\psi\p=1_S$).    We write $\wh{X^+}$ for the free
profinite semigroup generated by $X$.  %The following lemma is easily
%proved by a standard compactness argument.

\begin{Lemma}\label{projectivesemi}
Let $S$ be an $X$-generated finite projective semigroup.  Then the
canonical continuous projection
$\p:\wh{X^+}\to S$ splits and so $S$ embeds in $\wh{X^+}$.
\end{Lemma}
\begin{proof}
Write $\wh{X^+} = \ilim S_i$ where the $S_i$ are $X$-generated
finite semigroups;  we may assume that each $S_i$ maps onto $S$ via
the projection $\p_i:S_i\twoheadrightarrow S$.  It's easy to see
that the isomorphism $\mathrm{Hom}(S,\wh{X^+})\cong \ilim
\mathrm{Hom}(S,S_i)$ sends the splittings of $\p$ to the inverse
limit of the splittings of the $\p_i$.  As the inverse limit of
finite non-empty sets is non-empty, this completes the proof.
%Let $C_i$ be the set of splittings of $\p_i$. The $C_i$ form an
%inverse system of non-empty finite sets, whose non-empty inverse
%limit is the set of splittings of $\p$.
%Let $C_i$ be the set of all homomorphisms $S\to S_i$ splitting
%$\p_i$. By assumption $C_i\neq\emptyset$.  If $\pi_{ij}:S_i\to S_j$ denotes
%the canonical map and $f\in C_i$,
%then $f\pi_{ij}\in C_j$.  So the $C_i$ form an inverse system of
%non-empty sets.  Hence $\ilim C_i\neq \emptyset$.  It is easy to see
%that $\ilim C_i$ is just the set of homomorphisms $S\to \wh{X^+}$
%splitting $\p$.  Thus there is a homomorphism splitting $\p$.
\end{proof}

Since $\wh{X^+}$ embeds in the free profinite monoid on $X$, we
obtain as a consequence of Lemma~\ref{projectivesemi} and
Corollary~\ref{mainresult} the following theorem.

\begin{Thm}
Every finite projective semigroup is a band.
\end{Thm}

There are many examples of non-trivial projective
semigroups, such as chains of idempotents or left/right zero
semigroups~\cite{qtheor}.

\bibliographystyle{abbrv}
\bibliography{standard}

\end{document}